\documentclass[12pt,a4paper]{article}%
\usepackage{graphicx}
\usepackage{amsmath}%
\usepackage{amsfonts}%
\usepackage{amssymb}
\newtheorem{theorem}{Theorem}

\newtheorem{conjecture}[theorem]{Conjecture}

\newtheorem{lemma}[theorem]{Lemma}

\newenvironment{proof}[1][Proof]{\textbf{#1.} }{\ \rule{0.5em}{0.5em}}

\newcommand{\RR}{\ensuremath{\mathbb{R}}}
\newcommand{\ZZ}{\ensuremath{\mathbb{Z}}}
\newcommand{\ra}{\ensuremath{\rightarrow}}
\renewcommand{\P}{\ensuremath{\mathbb{P}}}
\newcommand{\E}{\ensuremath{\mathbb{E}}}
\DeclareMathOperator{\Var}{Var}
\DeclareMathOperator{\Cov}{Cov}
\newcommand{\qedsymbol}{\ensuremath{\blacksquare}}
\newcommand\ceil[1]{\lceil#1\rceil}
\newcommand{\ol}{\ensuremath{\overline}}

\begin{document}

\title{The diameter of a random Cayley graph of $Z_q$}

\author{{Gideon Amir \thanks{Department of Mathematics, University of Toronto,
Toronto ON, M5S 2E4, Canada. gidi.amir@gmail.com}} \quad {Ori Gurel-Gurevich
\thanks{Microsoft Research, One Microsoft Way, Redmond, WA 98052-6399, USA. origurel@microsoft.com}} }

\maketitle

\begin{abstract}
    Consider the Cayley graph of the cyclic group of prime order
    $q$ with $k$ uniformly chosen generators. For fixed $k$, we prove that the
    diameter of said graph is asymptotically (in $q$) of order
    $\sqrt[k]{q}$.

    The same also holds when the generating set is taken to be a symmetric set of size $2k$.
\end{abstract}

\section{Introduction}

Let $G$ be a finite group. Let $S$ be a subset of $G$. The
(directed) \emph{Cayley graph} of $G$ (w.r.t. $S$) is a graph
$(V,E)$ with $V=G$ and $(g,h)\in E$ if and only if $h^{-1}g \in S$.
The elements of $S$ are then called \emph{generators} and
$S$ the \emph{generating set}. If $S$ is symmetric (w.r.t. inversion)
then the resulting Cayley graph is essentially undirected - if
$(g,h)$ is in $E$ then so is $(h,g)$.

A "random random walk" on a group $G$ is a random walk on the
Cayley graph of $G$, with a generating set chosen randomly in some
fashion. The random walk itself may be simple, each edge having
equal probability at each step, or not simple, with some
nonuniform distribution on the edges. Usually, the generating set
is chosen uniformly from all sets of some prefixed size $k$.

Various aspects, most notably the typical mixing time, of random random
walks on different finite groups have been studied. (See
\cite{AlonRoichman,Wilson} for some examples, and
\cite{survey} which gives a comprehensive survey). The results
usually refer to Abelian groups in varying degrees of generality,
from cyclic groups up to general finite Abelian groups.

Roichman (\cite{Roichman}) notes that the diameter of the random
Cayley graph is bounded by a constant times the mixing time, and
applies this bound to the case of general groups of order $n$ and
$k=\lfloor\log^a n \rfloor$. The resulting bound is $\frac{a}{a-1}
\log_k n$, which is proved to be tight for the case of Abelian
groups.

For the cyclic group $\ZZ_n$, Hildebrand \cite{hildebrand} proved
that the mixing time is of order $n^{2/(k-1)}$, and therefore, this
is also a bound on the diameter of this random Cayley graph.
However, in contrast with the results in \cite{Roichman}, in this
case the diameter is actually much smaller - we prove it to be $O(n^{1/k})$.

Note that as far as mixing times are concerned, there is no
difference between a particular set of generators, $S$ and the set
$S+c$, attained by adding a constant $c\in \ZZ_q$ to all the
generators in $S$. The diameter, however, might change
significantly. To see this, consider, for example, a generating set
of two elements, $S=\{1,\lceil\sqrt{q}\rceil\}$ in $\ZZ_q$ (where
$q$ is prime). The diameter of this Cayley graph is $2\sqrt{q}$.
Now, observe $S'=S-1=\{0,\lfloor\sqrt{q}\rfloor\}$. The diameter of
this Cayley graph is now $q$. Put another way, the diameter does not
change when adding or removing $0$ from the generating set but the
mixing time might change considerably.

Another point of notice is the question of symmetric vs. asymmetric
generating sets. The results in \cite{hildebrand} are for asymmetric
generating sets, i.e. $S$ contains just the $k$ randomly chosen
generators. We might as well ask about the mixing times and
diameters w.r.t. $\ol{S}=S\bigcup (-S)$. The resulting random
walk is now symmetric, which is sometimes more natural to consider.
It seems that the results in \cite{hildebrand}, when applied to the
symmetric case, would yield a mixing time of order $n^{2/k}$. In
contrast, the results in this paper apply equally to the symmetric
case.

It should be noted that the asymptotic behavior of both the mixing
time and the diameter, both in the symmetric and asymmetric case
are the same as in the case of a $k$ dimensional tori of volume
$q$. This is not coincidental, the structure of the Cayley graph
of $G$ w.r.t. $S$ is actually that of $\ZZ^k$, modulo some
$k$-dimensional lattice, which contains the lattice of all
multiples of $q$. Perhaps the mixing time results of
\cite{hildebrand} could be proved in a more elementary manner
using that perspective.

\section{Main results and open questions}

Let $\ZZ_q$ be the cyclic group of order $q$, a prime number. Let
$g_1,..,g_k$ be $k$ random generators chosen uniformly and independently
from $\ZZ_q$. Denote by $Diam(q,k)$ the random variable which is the
diameter of the resulting (directed) Cayley graph. The same proofs work,
\textit{mutatis mutandis}, for the diameter of the \emph{undirected} Cayley graph,
that is, if our generating set is taken to be $\{g_1, -g_1, \ldots , g_k, -g_k\}$.

A simple counting argument shows that the diameter is at least $\Omega(\sqrt[k]{q})$. We
prove that the diameter is $\Theta(\sqrt[k]{q})$ in the following sense:

\begin{theorem}\label{ThUB}
For all $k>0$,
$$\lim_{C\ra\infty} \limsup_{q\ra\infty} \P(Diam(q,k) > C \sqrt[k]{q}) = 0$$
\end{theorem}

Also, this result is tight in the sense that:

\begin{theorem}\label{ThLB}
For all $k>0$ and all $C$,
$$\liminf_{q\ra\infty} \P(Diam(q,k) > C \sqrt[k]{q}) > 0$$
\end{theorem}

In other words, the limit behavior of the distribution of $\frac{Diam(q,k)}{\sqrt[k]{q}}$ is non-degenerate.
This seems to hint at the following conjecture:

\begin{conjecture}
$\frac{Diam(q,k)}{\sqrt[k]{q}}$ converges (in distribution) to some
distribution $D(k)$ on $\RR$ which has a non-compact support.
\end{conjecture}

If this conjecture is true, an obvious question would be to find out what this
limit distribution is.

\section{Proof of Theorem \ref{ThUB}}

For $x \in \ZZ_q$ and $\ol{i}=\{i_1,..,i_k\} \in \{0,\ldots,L\}^k$ a
vector of indices, let $A_{\ol{i}}^x$ be the event $i_1
g_1+i_2 g_2+..+i_k g_k = x \pmod{q}$. Let
$A^x_L=\bigcup_{\ol{i}\in L^k} A_{\ol{i}}^x$ and let
$A_L=\bigcap_{x \in Z_q} A^x_L$. We abuse the notation and identify an event with its indicator function.

If $A_L$ occurs then the diameter of the Cayley graph is at most $kL$, while if $A_L$ doesn't occur then the diameter is at least $L$, which is the same order of magnitude, since $k$ is fixed. Therefore, to prove both theorems it is enough to bound $\P(A_{C \sqrt[k]{q}})$ from above and below.

Obviously, for any $\ol{i}\neq 0^k$ and any $x\in \ZZ_q$ we have
$\E(A_{\ol{i}}^x)=1/q$.

Next we want to calculate $\E(A_{\ol{i}}^x
A_{\ol{j}}^x)$. This is the same as asking how many solutions, in $(\ZZ_q)^k$, are there for:
$$i_1 g_1+i_2 g_2+..+i_k g_k = x$$
$$j_1 g_1+j_2 g_2+..+j_k g_k = x$$
If $\ol{i}$ and $\ol{j}$ are linearly independent over
$\ZZ_q$ then the number of solutions is exactly $q^{k-2}$. In that
case $\E(A_{\ol{i}}^x A_{\ol{j}}^x)=1/q^2$ and
therefore the events are independent. If $\ol{i}$ and
$\ol{j}$ are linearly dependent over $\ZZ_q$ then
$\ol{i}=\lambda\ol{j}$ for some $\lambda\neq 1$. In
that case there are no solutions since $x\neq \lambda x$ (except
for $x=0$ which we can ignore). Therefore, in that case
$$\Cov(A_{\ol{i}}^x,
A_{\ol{j}}^x)=\E(A_{\ol{i}}^x
A_{\ol{j}}^x)-\E(A_{\ol{i}}^x)\E(A_{\ol{j}}^x)=0-1/q^2
< 0 \ .$$

Let $B^x_L=\sum_{\ol{i}\in L^k} A_{\ol{i}}^x$. Notice
that $A^x_L=0$ if and only if $B^x_L=0$. By linearity of expectation,
$$\E(B^x_L)=\sum_{\ol{i}\in L^k}
\E(A_{\ol{i}}^x)=\frac{L^k}{q} \ .$$

$\Var(A_i^x) = \frac{1}{q}(1-\frac{1}{q}) < \frac{1}{q}$ and
$\Cov(A_i^x,A_j^x) \leq 0$, giving
$$\Var(B^x)=\sum_{\ol{i}\in L^k} \Var(A_{\ol{i}}^x) + \sum_{\ol{i}\in L^k}
\sum_{\ol{i}\neq\ol{j}\in L^k}
\Cov(A_{\ol{i}}^x,A_{\ol{j}}^x) < \frac{L^k}{q} \ .$$

Chebyshev's inequality now yields
$$\P(B^x_L=0)\le \P(|B^x_L-\E(B^x_L)|\ge \E(B^x_L)) \le
\frac{\Var(B^x_L)}{\E(B^x_L)^2}<\frac{L^k/q}{(L^k/q)^2}=\frac{q}{L^k}$$
and therefore
$$\P(A^x_L)=1-\P(B^x_L=0)\geq 1-\frac{q}{L^k}$$

Let $T_L=\{x| A^x_L \}$ be the set of all points in $\ZZ_q$ that can be reached by using each generator at most $L$ times, and let $B_L=|T_L|=\sum_{x\in Z_q} A^x_L$.
Fix $C>0$ and let $L=C\sqrt[k]{q}$. We now have $\P(A^x_L)\geq 1-\frac{q}{L^k}=1-\frac{1}{C^k}$.  Therefore
$\E(B_L)\geq q(1-\frac{1}{C^k})$. Since $B\le q$, we can use Markov's inequality on $q-B_L$ to get
$$\P(B_L>\frac{q}{2})=1-\P(q-B_L > \frac{q}{2})\geq 1- \frac{\frac{q}{C^k}}{\frac{q}{2}}=1-\frac{2}{C^k} \ .$$

Now if $B_L>\frac{q}{2}$ then for every $x\in \ZZ_q$ we have $T\cap(x-T)\neq \emptyset$. This means that $A_{2L}$ occurs and therefore the diameter is at most $2kL$.

Therefore,
$$\P(Diam(q,k) > C \sqrt[k]{q})\le \P(B_{(C/2k) \sqrt[k]{q} } > \frac{q}{2}) \le \frac{2}{(C/2k)^k} \underset{C\ra \infty}{\longrightarrow} 0$$
as required. \qedsymbol

\section{Proof of Theorem \ref{ThLB}}

Fix some $D<1$ and let $L= D\sqrt[k]{q}$. Consider the events $A_{\ol{i}}^0$ and
$A^0_L$ as previously defined. As before, if $\ol{i}$ and
$\ol{j}$ are linearly independent over $Z_q$ then
$A_{\ol{i}}^0$ and $A_{\ol{j}}^0$ are independent
events. If $\ol{i}$ and $\ol{j}$ are linearly
dependent then $A_{\ol{i}}^0$ and $A_{\ol{j}}^0$ are
in fact the same event. How many distinct events do we have among
$\{A_{\ol{i}}^0\}_{i \in L^k}$?

\begin{lemma}
There are at least $L^k/2=D^k q/2$ such distinct events.
\end{lemma}

\begin{proof}
Let $\ol{i}$ and $\ol{j}$ be linearly dependent, i.e.
there exist $c\in Z_q$ such that $\ol{i}=c\ol{j}$. In
particular $i_0=c j_0 \pmod{q}$ and $i_1=c j_1 \pmod{q}$.
eliminating $c$, we get $i_0 j_1 = i_1 j_0 \pmod{q}$. Since
$i_0$,$i_1$,$j_0$ and $j_1$ are all less than $\sqrt{q}$ (since $k\ge2$) we get
$i_0 j_1 = i_1 j_0$.

Take all $\ol{i} \in L^k$ for which $i_0$ and $i_1$ are
coprime in $\ZZ$. If $i_0$ and $i_1$ are coprime and
$j_0$ and $j_1$ are coprime and $i_0 j_1 = j_0 i_1$ then $i_0=j_0$
and $i_1=j_1$. Therefore, among the vectors considered above every two
are linearly independent, so the corresponding events are distinct.

Given $L$ how many pairs $i_0,i_1 <L$ are coprime? It is a well
known fact (see \cite{penguin}) that the fraction of coprime pairs tends to, and is always
greater than, $6/\pi^2 > 1/2$.
\end{proof}

Let $I \subset L^k$ be a set of index vectors such that every two
are linearly independent and $|I|= \ceil{D^kq/2}L$. Let
$X=\sum_{\ol{i}\in I} A_{\ol{i}}^0$ and notice that
$\P(A^0_L)\geq \P(X>0)$.

\begin{lemma}\label{D3}
$\P(X>0)>D^k/3$
\end{lemma}
\begin{proof}
For all $\ol{i}\in I$ we have $\E(A_{\ol{i}}^0)=\frac{1}{q}$ so
$$\E(X)=\frac{D^k q}{2} \frac{1}{q}=\frac{D^k}{2}$$
and these events are pairwise independent, so
$$\E(X^2)=\frac{D^k q}{2} \frac{1}{q} + \frac{D^k q}{2}(\frac{D^k q}{2}-1) \frac{1}{q^2}=\frac{D^k}{2}(1-\frac{1}{q}+\frac{D^k}{2})$$

Since $X$ is nonnegative, from Cauchy-Schwartz inequality we get
$$\E(X)^2=\E(X 1_{X>0})^2 \leq \E(X^2)\E((1_{X>0})^2)=\E(X^2)\P(X>0)$$
and therefore
$$\P(X>0)\geq \frac{\E(X)^2}{\E(X^2)}=\frac{\frac{D^{2k}}{4}}{\frac{D^k}{2}(1-\frac{1}{q}+\frac{D^k}{2})}
\geq \frac{D^k}{3}$$
as required.
\end{proof}

From lemma \ref{D3} we get that the probability of $A^0_L$ is
bounded away from 0 regardless of $q$. Next we shall show that if
$A^0_L$ occurs then many different $\ol{i}$ yield the same
member of $Z_q$, in which case the diameter cannot be too small.

\begin{lemma}\label{LB}
Let $C$ be such that $k D C^{k-1} < 1$. If $A^0_{D \sqrt[k]{q}}$ occurs then
$$Diam(q,k) > C\sqrt[k]{q} \ .$$
\end{lemma}
\begin{proof}
Let $L=D\sqrt[k]{q}$ and let $\ol{i} \in \{0,\ldots,L\}^k$ be such that
$A_{\ol{i}}^0$ occurs. If $\ol{j}$ and $\ol{j'}$
differ by a multiple of $\ol{i}$ then
$$j_1 g_1+j_2 g_2+..+j_k g_k=j'_1 g_1+j'_2 g_2+..+j'_k g_k \pmod{q} \ .$$

Therefore for every $\ol{j} \in \{0,\ldots,L\}^k$ there
exists $\ol{j'}$ such that
$$j_1 g_1+j_2 g_2+..+j_k g_k=j'_1 g_1+j'_2 g_2+..+j'_k g_k \pmod{q}$$
and $j'_n \leq D\sqrt[k]{q}$ for some $1\le n \leq k$. The number
of such $\ol{j'}$ is bounded by $kD\sqrt[k]{q}(C\sqrt[k]{q})^{k-1}=k D C^{k-1} q < q$.
Therefore, if $A^0_{D\sqrt[k]{q}}$ occurs not all vertices are covered by combinations in $\{0,\ldots,C\sqrt[k]{q}^k$ and hence the diameter is at least $C \sqrt[k]{q}$.
\end{proof}

To wrap up the proof, given $C$, let $D = 1/(2kC^{k-1})$. From lemma \ref{D3} we conclude that $A^0_{D\sqrt[k]{q}}$ occurs with probability at least $D/3$, in which
case, by lemma \ref{LB} we get $Diam(q,k) > \sqrt[k]{Cq}$. \qedsymbol

\subsection*{Acknowledgements}

The authors thank Itai Benjamini, for suggesting this problem and for useful discussions.


\begin{thebibliography}{99}

\bibitem{AlonRoichman}
N. Alon, Y. Roichman (1994) Random Cayley graphs and expanders.
\emph{Random Structures and Algorithms}, \textbf{5}(2), 271-284

\bibitem{survey}
M. V. Hildebrand (2005) A survey of random random walks on finite
groups. \emph{Probability Surveys 2}, 33-63

\bibitem{hildebrand}
M. V. Hildebrand (1994) Random walks supported on random points of
$\ZZ/n\ZZ$. \emph{Probability Theory and Related Fields} 100,
191-203


\bibitem{penguin}
D. Wells (1986) The Penguin Dictionary of Curious and Interesting
Numbers. Middlesex, England: Penguin Books, 28-29

\bibitem{Wilson}
D.B. Wilson, Random random walks on $\ZZ^d_2$. \emph{Probability
Theory and Related Fields} 108 (1997), no. 4

\bibitem{Roichman}
Y. Roichman, On random random walks. \emph{Ann. Probab.} 24
(1996), no. 2, 1001--1011.




\end{thebibliography}
\end{document}